\documentclass{amsart}
\newtheorem{theorem}{Theorem}[section]
\newtheorem{lemma}[theorem]{Lemma}

\theoremstyle{definition}
\newtheorem{definition}[theorem]{Definition}
\newtheorem{example}[theorem]{Example}

\theoremstyle{remark}
\newtheorem{remark}[theorem]{Remark}

\numberwithin{equation}{section}

\begin{document}
 \title{Rohit self-class group}
\author{Author One}
\curraddr{\textbf{Rohit Madukar Patne} M.Sc. mathematics and computing ,IIT Guwahati,Guwahati,India-781039}
\curraddr{Home address:\textbf{Rohit Madhukar Patne},plot no 51,Ramkrushna Nagar,behind raghute bhavan,Dighori,\\
Umrer Road ,Nagpur\\Maharashtra\\India-440034}
\email{m.patne@iitg.ernet.in or m.patne@yahoo.com}
\thanks{In advance,I would like to thank all respected  sir who will read this article,and also for giving me helpfull  suggestion}
\date{27 nov 2011}
\dedicatory{This paper is dedicated to my mother and father,litle brother sushil,all my teacher,and my sir A.S.Muktibhodh }
\begin{abstract}In this paper, I will define \textbf{self-class} and \textbf{conjugacy self-class}.and also give
 some result about it,then i will define \textbf{self-class group}.This paper is totally illustrative means 
I have tried to explain each definition and theorem 
by taking sutaible example.If anyone will read this,then he will observe 
the beauty of \textbf{self-class} and \textbf{conjugacy self-class}
in study of Group theory module theory and field theory and other brances of mathematics .
and I also discuss \textbf{non-selfclass},and I have given the conjecture about non-self class Group. 
% aim of this article is to study group property 
%how nice structure will obtain
 \end{abstract}
\maketitle
%%%%%%%%%%%%%%%%%%%%%%%%%%%%%%%%%%%%%%%%%%%%%%%%%%%%%%%%%%%%%%%%%%%%%%%%%%%%%%%%%%%%%%%%%%%%%%%%%%%%%%%%%%%%%%%%%%%%%%%%%%%%%%%%%%%%%%%%%%%%%%%%%%%%%%%%%%%%%%%%
\section{Introduction}
Let $G$ be a finite group, and $K$ a abelian proper normal subgroup of $G$. Then there
are two natural partitions of $G$. One partition is given by the conjugacy
classes, while the other one by the cosets of $K$. Here \textbf{self-class} and \textbf{conjugacy self-class} are another two way of 
partition group G.also one can see that\textbf{self-class} is different from Double cosets.
Detail about \textbf{self-class} and \textbf{conjugacy self-class} are discuss in this paper.

%%%%%%%%%%%%%%%%%%%%%%%%%%%%%%%%%%%%%%%%%%%%%%%%%%%%%%%%%%%%%%%%%%%%%%%%%%%%%%%%%%%%%%%%%%%%%%%%%%%%%%%%%%%%%%%%%%%%%%%%%%%%%%%%%%%%%%%%%%%%%%%%%%%%%%%%%%%%%
\section{some definition}
\begin{remark}
 I have assume that \textbf{ \textasteriskcentered} is operation on Group $G$ ,hence instead of writting $(G$,\textasteriskcentered$)$ I will use simply 
Group $G$. if other operation 
will come I will mention that operation .
\end{remark}
\begin{definition} Let $G$ be a group and  $H$ be a abelian sugroup of $G$ ,
Then the relation of the form $b$=$x$\textasteriskcentered$a$\textasteriskcentered$x$
  such that $x$ $\in$ $H$ and $a$,$b$ $\in$ $G$ 
 %\\$a$\textasteriskcentered$x$=$x$\textasteriskcentered$a$ $\forall$ $x$ $\in$ $H$ and $a$ $\in$ $G$ 
is an equivalence relation and it is denoted by $b$ $\sim$ $a$.
        \end{definition}
\begin{proof}
 Let $a$,$b$,$c$ $\in$ $G$ ,and $H$ be a abelian sugroup of $G$ such that $\exists$ $x$ $\in$ $H$ such that
$b$=$x$\textasteriskcentered$a$\textasteriskcentered$x$	\\
$\Rightarrow$ $b$ $\sim$ $a$\\
also $\exists$ $y$ $\in$ $H$ such that\\
$a$=$y$\textasteriskcentered$c$\textasteriskcentered$y$,\\
$\Rightarrow$ $a$ $\sim$ $c$\\
clearly,
$b$=$x$\textasteriskcentered$y$\textasteriskcentered$c$\textasteriskcentered$x$\textasteriskcentered$y$,since $H$ be a abelian sugroup of $G$ \\
$\Rightarrow$ $b$ $\sim$ $c$\\
from this we can say that, $G$ be a group and  $H$ be a abelian sugroup of $G$ ,
Then the relation of the form $b$=$x$\textasteriskcentered$a$\textasteriskcentered$x$
  such that $x$ $\in$ $H$ and $a$,$b$ $\in$ $G$ is an equivalence relation. 	
\end{proof}

\begin{definition}self-class is defined as\\ 
$g$ $($ $a$ $)$=$\{$ $x$\textasteriskcentered$a$\textasteriskcentered$x$ $/$  $\forall$ $x$ $\in$ $H$,$a$ $\in$ $G$,$H$ 
be a abelian sugroup of $G$ ,also Here \textbf{ \textasteriskcentered} is an operation which is equivalent to operation define on group $G$. %and $a$\textasteriskcentered$x$=$x$\textasteriskcentered$a$ $\forall$ $x$ $\in$ $H$ and $a$ $\in$ $G$ 
$\}$ 
\end{definition}
\begin{remark}
 Notation use to denote self-class of $a$ $\in$ $G$ is $g$ $($ $a$ $)$.
\end{remark}
 
%\begin{remark}
 %I have assume that \textbf{ \textasteriskcentered} is operation on Group $G$ ,hence instead of writing $(G$,\textasteriskcentered$)$ I will use simply 
%Group $G$. if other operation 
%will come I will mention that operation .
%\end{remark}

%\begin{remark} 
 %Here \textbf{ \textasteriskcentered} is an operation which is equivalent to operation define on group $G$. I will use this operation ,if other operation 
%will come I will mention that operation .
%\end{remark}
\begin{example} \textbf{To illustrate the above definition by using example}\\ Let $G$ be a abelian group of oder $6$ such that\\ $G$=$\{$ $e,$ $a$,$a^{2}$,$a^{3}$,$a^{4}$,$a^{5}$ $\}$
\\ then clearly $H$=$\{$ $e,$ $a^{3}$ $\}$ and $K$=$\{$ $e,$ $a^{2}$,$a^{4}$ $\}$ are proper ablelian subgroup of group $G$.
Then self-class  of Group $G$ over $H$ is as follow,\\
 $g$ $($ $a$ $)$=$\{$ $e$\textasteriskcentered$a$\textasteriskcentered$e$,$a^{3}$\textasteriskcentered$a$\textasteriskcentered$a^{3}$ $\}$\\
$\Rightarrow$ $g$ $($ $a$ $)$=$\{$ $a$ $\}$.\\
 $g$ $($ $e$ $)$=$\{$ $e$\textasteriskcentered$e$\textasteriskcentered$e$,$a^{3}$\textasteriskcentered$e$\textasteriskcentered$a^{3}$ $\}$\\
$\Rightarrow$ $g$ $($ $e$ $)$=$\{$ $e$ $\}$.\\
 $g$ $($ $a^{2}$ $)$=$\{$ $e$\textasteriskcentered$a^{2}$\textasteriskcentered$e$,$a^{3}$\textasteriskcentered$a^{2}$\textasteriskcentered$a^{3}$ $\}$\\
$\Rightarrow$ $g$ $($ $a^{2}$ $)$=$\{$ $a^{2}$ $\}$.\\
$g$ $($ $a^{3}$ $)$=$\{$ $e$\textasteriskcentered$a^{3}$\textasteriskcentered$e$,$a^{3}$\textasteriskcentered$a^{3}$\textasteriskcentered$a^{3}$ $\}$\\
$\Rightarrow$ $g$ $($ $a^{3}$ $)$=$\{$ $a^{3}$ $\}$.\\
$g$ $($ $a^{4}$ $)$=$\{$ $e$\textasteriskcentered$a^{4}$\textasteriskcentered$e$,$a^{3}$\textasteriskcentered$a^{4}$\textasteriskcentered$a^{3}$ $\}$\\
$\Rightarrow$ $g$ $($ $a^{4}$ $)$=$\{$ $a^{4}$ $\}$.\\
 $g$ $($ $a^{5}$ $)$=$\{$ $e$\textasteriskcentered$a^{5}$\textasteriskcentered$e$,$a^{3}$\textasteriskcentered$a^{5}$\textasteriskcentered$a^{3}$ $\}$\\
$\Rightarrow$ $g$ $($ $a^{5}$ $)$=$\{$ $a^{5}$ $\}$.\\
hence the collection of all self-class of Group $G$ over $H$ is set $S_{H}$ such that\\
 $S_{H}$ =$\{$  $g$ $($ $e$ $)$, $g$ $($ $a$ $)$, $g$ $($ $a^{2}$ $)$, $g$ $($ $a^{3}$ $)$, $g$ $($ $a^{4}$ $)$, $g$ $($ $a^{5}$ $)$ $\}$\\
 \textbf{also} we can also find the self-class of Group $G$ over $K$ as follow\\%%%%%%%%%%%%%%%%%%%%%%%%%%%%%%%%%%%%%%%%%%%%%%%%%%%%%%%%%%%%%
 $g$ $($ $a$ $)$=$\{$ $e$\textasteriskcentered$a$\textasteriskcentered$e$,$a^{2}$\textasteriskcentered$a$\textasteriskcentered$a^{2}$,
         $a^{4}$\textasteriskcentered$a$\textasteriskcentered$a^{4}$ $\}$\\
$\Rightarrow$ $g$ $($ $a$ $)$=$\{$ $a$,$a^{5}$,$a^{3}$ $\}$.\\
 $g$ $($ $a^{2}$ $)$=$\{$ $e$\textasteriskcentered$a^{2}$\textasteriskcentered$e$,$a^{2}$\textasteriskcentered$a^{2}$\textasteriskcentered$a^{2}$,
         $a^{4}$\textasteriskcentered$a^{2}$\textasteriskcentered$a^{4}$ $\}$\\
$\Rightarrow$ $g$ $($ $a$ $)$=$\{$ $a$,$a^{2}$,$a^{4}$ $\}$.\\
clearly, one can obeserve that  $g$ $($ $a$ $)$=  $g$ $($ $a^{3}$ $)$= $g$ $($ $a^{5}$ $)$ and\\
$g$ $($ $e$ $)$=  $g$ $($ $a^{2}$ $)$= $g$ $($ $a^{4}$ $)$\\
hence the collection of all self-class of Group $G$ over $K$ is set $S_{K}$ such that\\
$S_{K}$ =$\{$ $g$ $($ $a$ $)$,$g$ $($ $e$ $)$ $\}$.
\end{example}

%\begin{lemma}
 %          collection  of  all $g$ $($ $a$ $)$  form a group  $($say $K$ $)$
 %with operation\\ $g$ $($ $a$ $)$\textasteriskcentered$g$ $($ $b$ $)$=$\{$ c$.$d$/$ $\forall$ $c$ $\in$ $g$ $($ $a$ $)$ and $\forall$ $d$ $\in$ $g$ $($ $b$ $)$ and dot means $.$
%is an operation which is equal to operation on $G$ and \textasteriskcentered is an operation on $K$ $\}$
%\\iff $a$ $ .$ $x$=$x$ $.$ $a$ $\forall$ $x$ $\in$ $H$ and $a$ $\in$ $G.$
 %         \end{lemma}
%%%%%%%%%%%%%%%%%%%%%%%%%%%%%%%%%%%%%%%%%%%%%%%%%%%%%%%%%%%%%%%%%%%%%%%%%%%%%%%%%%%%%%%%%%%%%%%%%%%%%%%%%%%%%%%%%%%%%%%%%%%%%%%%%%%%%%%%%%%%%%%%%%%%%%%%%%%%%%%%%%%%%%
\begin{definition} Let $G$ be a group and  $H$ be a sugroup of $G$ ,
Then the relation of the form $b$=$x^{-1}$\textasteriskcentered$a$\textasteriskcentered$x$
  such that $x$ $\in$ $H$ and $a$,$b$ $\in$ $G$ and it is denoted by $b$ $\sim$ $a$.
 %\\$a$\textasteriskcentered$x$=$x$\textasteriskcentered$a$ $\forall$ $x$ $\in$ $H$ and $a$ $\in$ $G$ 
is an equivalence relation.
        \end{definition}
\begin{proof}
 one can easily check that the given relation is an equivalence relation on group $G$.
\end{proof}

\begin{definition}conjugacy-selfclass is defined as\\ 
$cl_{s_{H}}$ $($ $a$ $)$=$\{$ $x^{-1}$\textasteriskcentered$a$\textasteriskcentered$x$ $/$  $\forall$ $x$ $\in$ $H$,$x^{-1}$ $\in$ $H$ ,$a$ $\in$ $G$,$H$ 
be a sugroup of $G$,and also \textbf{ \textasteriskcentered} is an operation which is equivalent to operation define on group $G$. %and $a$\textasteriskcentered$x$=$x$\textasteriskcentered$a$ $\forall$ $x$ $\in$ $H$ and $a$ $\in$ $G$ 
$\}$ 
\end{definition}
\begin{remark}
 Notation use to denote conjugacy-selfclass of $a$ $\in$ $G$ is $cl_{s_{H}}$ $($ $a$ $)$.
\end{remark}
\begin{example}\textbf{To illustrate the above definition by using example}\\
Let $G$ be a abelian group of oder $6$ such that\\ $G$=$\{$ $e,$ $a$,$a^{2}$,$a^{3}$,$a^{4}$,$a^{5}$ $\}$
\\ then clearly $H$=$\{$ $e,$ $a^{3}$ $\}$ and $K$=$\{$ $e,$ $a^{2}$,$a^{4}$ $\}$ are subgroup of group $G$.
Then conjugacy-selfclass of Group $G$ over $H$ is as follow,\\
 $cl_{s_{H}}$ $($ $a$ $)$=$\{$ $e$\textasteriskcentered$a$\textasteriskcentered$e$,$a^{3}$\textasteriskcentered$a$\textasteriskcentered$a^{3}$ $\}$\\
$\Rightarrow$ $cl_{s_{H}}$ $($ $a$ $)$=$\{$ $a$ $\}$.\\
 $cl_{s_{H}}$ $($ $e$ $)$=$\{$ $e$\textasteriskcentered$e$\textasteriskcentered$e$,$a^{3}$\textasteriskcentered$e$\textasteriskcentered$a^{3}$ $\}$\\
$\Rightarrow$ $cl_{s_{H}}$ $($ $e$ $)$=$\{$ $e$ $\}$.\\
 $cl_{s_{H}}$ $($ $a^{2}$ $)$=$\{$ $e$\textasteriskcentered$a^{2}$\textasteriskcentered$e$,$a^{3}$\textasteriskcentered$a^{2}$\textasteriskcentered$a^{3}$ $\}$\\
$\Rightarrow$ $cl_{s_{H}}$ $($ $a^{2}$ $)$=$\{$ $a^{2}$ $\}$.\\
$cl_{s_{H}}$ $($ $a^{3}$ $)$=$\{$ $e$\textasteriskcentered$a^{3}$\textasteriskcentered$e$,$a^{3}$\textasteriskcentered$a^{3}$\textasteriskcentered$a^{3}$ $\}$\\
$\Rightarrow$ $cl_{s_{H}}$ $($ $a^{3}$ $)$=$\{$ $a^{3}$ $\}$.\\
$cl_{s_{H}}$ $($ $a^{4}$ $)$=$\{$ $e$\textasteriskcentered$a^{4}$\textasteriskcentered$e$,$a^{3}$\textasteriskcentered$a^{4}$\textasteriskcentered$a^{3}$ $\}$\\
$\Rightarrow$ $cl_{s_{H}}$ $($ $a^{4}$ $)$=$\{$ $a^{4}$ $\}$.\\
 $cl_{s_{H}}$ $($ $a^{5}$ $)$=$\{$ $e$\textasteriskcentered$a^{5}$\textasteriskcentered$e$,$a^{3}$\textasteriskcentered$a^{5}$\textasteriskcentered$a^{3}$ $\}$\\
$\Rightarrow$ $cl_{s_{H}}$ $($ $a^{5}$ $)$=$\{$ $a^{5}$ $\}$.\\
hence the collection of all conjugacy-selfclass  of Group $G$ over $H$ is set $S_{H}$ such that\\
 $S_{H}$ =$\{$  $cl_{s_{H}}$ $($ $e$ $)$, $cl_{s_{H}}$ $($ $a$ $)$, $cl_{s_{H}}$ $($ $a^{2}$ $)$, $cl_{s_{H}}$ $($ $a^{3}$ $)$, $cl_{s_{H}}$ $($ $a^{4}$ $)$,
 $cl_{s_{H}}$ $($ $a^{5}$ $)$ $\}$\\
\textbf{also} we can also find the conjugacy-selfclass of Group $G$ over $K$ as follow\\

$cl_{s_{H}}$ $($ $a$ $)$=$\{$ $e$\textasteriskcentered$a$\textasteriskcentered$e$,$a^{2}$\textasteriskcentered$a$\textasteriskcentered$a^{4}$,
         $a^{4}$\textasteriskcentered$a$\textasteriskcentered$a^{2}$ $\}$\\
$\Rightarrow$  $cl_{s_{H}}$ $($ $a$ $)$=$\{$ $a$ $\}$.\\

 $cl_{s_{H}}$ $($ $a^{2}$ $)$=$\{$ $e$\textasteriskcentered$a^{2}$\textasteriskcentered$e$,$a^{2}$\textasteriskcentered$a^{2}$\textasteriskcentered$a^{4}$,
         $a^{4}$\textasteriskcentered$a$\textasteriskcentered$a^{2}$ $\}$\\
$\Rightarrow$  $cl_{s_{H}}$ $($ $a$ $)$=$\{$ $a^{2}$ $\}$.\\

$cl_{s_{H}}$ $($ $a^{3}$ $)$=$\{$ $e$\textasteriskcentered$a^{3}$\textasteriskcentered$e$,$a^{2}$\textasteriskcentered$a^{3}$\textasteriskcentered$a^{4}$,
         $a^{4}$\textasteriskcentered$a^{3}$\textasteriskcentered$a^{2}$ $\}$\\
$\Rightarrow$  $cl_{s_{H}}$ $($ $a$ $)$=$\{$ $a^{3}$ $\}$.\\
$cl_{s_{H}}$ $($ $a^{4}$ $)$=$\{$ $e$\textasteriskcentered$a^{4}$\textasteriskcentered$e$,$a^{2}$\textasteriskcentered$a^{4}$\textasteriskcentered$a^{4}$,
         $a^{4}$\textasteriskcentered$a^{4}$\textasteriskcentered$a^{2}$ $\}$\\
$\Rightarrow$  $cl_{s_{H}}$ $($ $a^{4}$ $)$=$\{$ $a^{4}$ $\}$.\\

$cl_{s_{H}}$ $($ $a^{5}$ $)$=$\{$ $e$\textasteriskcentered$a^{5}$\textasteriskcentered$e$,$a^{2}$\textasteriskcentered$a^{5}$\textasteriskcentered$a^{4}$,
         $a^{4}$\textasteriskcentered$a^{5}$\textasteriskcentered$a^{2}$ $\}$\\
$\Rightarrow$  $cl_{s_{H}}$ $($ $a^{5}$ $)$=$\{$ $a^{5}$ $\}$.\\

$cl_{s_{H}}$ $($ $e$ $)$=$\{$ $e$\textasteriskcentered$e$\textasteriskcentered$e$,$a^{2}$\textasteriskcentered$e$\textasteriskcentered$a^{4}$,
         $a^{4}$\textasteriskcentered$e$\textasteriskcentered$a^{2}$ $\}$\\
$\Rightarrow$  $cl_{s_{H}}$ $($ $e$ $)$=$\{$ $e$ $\}$.\\

hence the collection of all conjugacy-selfclass  of Group $G$ over $K$ is set $S_{K}$ such that\\
 $S_{K}$ =$\{$  $cl_{s_{H}}$ $($ $e$ $)$, $cl_{s_{H}}$ $($ $a$ $)$, $cl_{s_{H}}$ $($ $a^{2}$ $)$, $cl_{s_{H}}$ $($ $a^{3}$ $)$, $cl_{s_{H}}$ $($ $a^{4}$ $)$,
 $cl_{s_{H}}$ $($ $a^{5}$ $)$ $\}$\\

\end{example}

\begin{definition}
  Let $G$ be a group and $H$ be an abelian subgroup of $G$, $H$ is an \textbf{faithfull} subgroup of $G$ ,if 
	  collection of all $g$ $($ $a$ $)$ form a group,otherwise $H$ is called \textbf{nonfaithfull} subgroup of $G.$
        \end{definition}
\begin{example}\textbf{example of \textbf{faithfull} subgroup of $G$ and \textbf{nonfaithfull} subgroup of $G.$}\\
Here I will choose the non-abelian Group of oder 6 for explaining above term,for this example ,I have use [3],reader can take any example,
Let $G$ be an non-abelian Group of oder 6 such that element of $G$ is as follow,\\
$G$=$\{$ $E$,$A$,$B$,$C$,$D$,$K$ $\}$ clearly, $H^{'}$=$\{$ $A$,$E$ $\}$,$H^{''}$=$\{$ $C$,$E$ $\}$,and $H^{''}$=$\{$ $K$,$E$ $\}$ 
$H^{'''}$=$\{$ $D$,$B$ ,$E$ $\}$ are the abelian subgroup of Group $G$.
But collection of all self-class which is define over above subgroup doesnot form a Group 
\\,hence $H^{'}$,$H^{''}$,$H^{'''}$ are \textbf{nonfaithfull} subgroup of $G.$\\
and also $H^{''''}$=$\{$ $E$ $\}$ is only the \textbf{faithfull} subgroup of $G$.
one can easily check by apply definition of self-class.
\end{example}

\begin{definition}%Let $S$ be a set of collection of all self-class of Group $G$ over $H$
%If collection of all $g$ $($ $a$ $)$ form a group 
Let collection of all self-class of Group $G$ over $H$ is set $S_{H}$ such that  $S_{H}$ form a Group with respect to operation \textasteriskcentered,
then $G$ is called as \textbf{Rohit self-classes group}
           \end{definition}
\begin{definition}
          If $G$ be a group,$\exists$ special abelian subgroup $H$ of $G$ suchthat collection of all self-classes define over $H$
          form a group $K_{H}$ which is  isomorphic to  a group $M_{G}$  suchthat $M_{G}$ is an collection of all self-classes define over $G$ ,
	  then group $G$ is said tobe \textbf{$R$-group}.   
         \end{definition}
\begin{example}\textbf{illustration of above definition with example}\\
Let $G$ be a abelian group of oder $6$ such that\\ $G$=$\{$ $e,$ $a$,$a^{2}$,$a^{3}$,$a^{4}$,$a^{5}$ $\}$
\\  $H$=$\{$ $e,$ $a^{2}$,$a^{4}$ $\}$ are proper ablelian subgroup of group $G$.\\

the self-class of Group $G$ over $H$ as follow\\%%%%%%%%%%%%%%%%%%%%%%%%%%%%%%%%%%%%%%%%%%%%%%%%%%%%%%%%%%%%%
 $g$ $($ $a$ $)$=$\{$ $e$\textasteriskcentered$a$\textasteriskcentered$e$,$a^{2}$\textasteriskcentered$a$\textasteriskcentered$a^{2}$,
         $a^{4}$\textasteriskcentered$a$\textasteriskcentered$a^{4}$ $\}$\\
$\Rightarrow$ $g$ $($ $a$ $)$=$\{$ $a$,$a^{5}$,$a^{3}$ $\}$.\\
 $g$ $($ $a^{2}$ $)$=$\{$ $e$\textasteriskcentered$a^{2}$\textasteriskcentered$e$,$a^{2}$\textasteriskcentered$a^{2}$\textasteriskcentered$a^{2}$,
         $a^{4}$\textasteriskcentered$a^{2}$\textasteriskcentered$a^{4}$ $\}$\\
$\Rightarrow$ $g$ $($ $a$ $)$=$\{$ $a$,$a^{2}$,$a^{4}$ $\}$.\\
clearly, one can obeserve that  $g$ $($ $a$ $)$=  $g$ $($ $a^{3}$ $)$= $g$ $($ $a^{5}$ $)$ and\\
$g$ $($ $e$ $)$=  $g$ $($ $a^{2}$ $)$= $g$ $($ $a^{4}$ $)$\\
Then  collection of all self-class of Group $G$ over $H$ is set $K_{H}$
suchthat\\$K_{H}$=$\{$ $g$ $($ $a$ $)$,$g$ $($ $e$ $)$ $\}$ \\clearly $K_{H}$ form a group with respect to operation \textbullet  suchthat
\textbullet  is define as follow\\
$g$ $($ $a$ $)$ \textbullet $g$ $($ $b$ $)$=$\{$ c\textasteriskcentered d$/$ $\forall$ $c$ $\in$ $g$ $($ $a$ $)$ and
$\forall$ $d$ $\in$ $g$ $($ $b$ $)$ and  \textasteriskcentered
is an operation which is equal to operation on $G$ and \textbullet is an operation on $K_{H}$ $\}$.\\%%%%%%%%%%%%%%%%%%%%%%%%%%%%%%%%%%%%%%%%%%%%%%%%%%%%%
\textbf{also one} can check that collection of all self-class of Group $G$ over $G$ is set $M_{G}$
suchtat \\$M_{G}$=$\{$ $g$ $($ $a$ $)$,$g$ $($ $e$ $)$ $\}$ \\clearly $M_{G}$ form a group with respect to operation \textbullet  suchthat
\textbullet  is define as follow\\
$g$ $($ $a$ $)$ \textbullet $g$ $($ $b$ $)$=$\{$ c\textasteriskcentered d$/$ $\forall$ $c$ $\in$ $g$ $($ $a$ $)$ and
$\forall$ $d$ $\in$ $g$ $($ $b$ $)$ and  \textasteriskcentered
is an operation which is equal to operation on $G$ and \textbullet is an operation on $M_{G}$ $\}$.\\
\textbf{clearly} since $K_{H}$ and $M_{G}$ are group of oder 2,as we know that $\exists$ unique group of oder 2.\\
$K_{H}$ and $M_{G}$ are isomorphic.\\
$\Rightarrow$ group $G$ is an \textbf{$R$-group}.

 %\begin{remark}
  %one can easily calculate all self-class of Group $G$ over $G$ .
 %\end{remark}

\end{example}
\begin{remark}
  one can easily calculate all self-class of Group $G$ over $G$ .
 \end{remark}
\begin{remark}
I hope that reader are now familiar with self-class and conjugacy-selfclass,he/she can take other example whatever he/she know and 
apply the definition of self-class and conjugacy-selfclass .In next section, I will describe some result associated with self class.
any one can check easily.
\end{remark}

\section{some lemma and theorem}
\begin{lemma}
 Let $G$ be a group and  $H$ be a abelian sugroup of $G$,then $g$ $($ $a$ $)$ is define over $H$,then 
1$\leq$ $|$ $g$ $($ $a$ $)$ $|$ $\leq$ $|$ $H$ $|$.where $|$ $g$ $($ $a$ $)$ $|$ is cardinality of $g$ $($ $a$ $)$ means number of element present inside $g$ $($ $a$ $)$.
\end{lemma}
\begin{lemma}
 Let $G$ be a group and  $H$ be a sugroup of $G$,then $cl_{s_{H}}$ $($ $a$ $)$ is define over $H$,then 
1$\leq$ $|$ $cl_{s_{H}}$ $($ $a$ $)$ $|$ $\leq$ $|$$H$$|$.where $|$ $cl_{s_{H}}$ $($ $a$ $)$ $|$ is cardinality of$cl_{s_{H}}$ $($ $a$ $)$
  means number of element present inside $cl_{s_{H}}$ $($ $a$ $)$.
\end{lemma}
\begin{lemma}Let $G$ be a group and  $H$ be a abelian sugroup of $G$,then $g$ $($ $a$ $)$ is define over $H$ suchthat  
           collection  of  all $g$ $($ $a$ $)$  form a group  $($say $K$ $)$
 with operation\\ $g$ $($ $a$ $)$\textbullet$g$ $($ $b$ $)$=$\{$ c\textasteriskcentered d $/$ $\forall$ $c$ $\in$ $g$ $($ $a$ $)$ and $\forall$ $d$ $\in$ $g$ $($ $b$ $)$ and dot means \textasteriskcentered
is an operation which is equal to operation on $G$ and \textbullet is an operation on $K$ $\}$
\\iff $a$ \textasteriskcentered $x$=$x$ \textasteriskcentered $a$ $\forall$ $x$ $\in$ $H$ and $a$ $\in$ $G.$
          \end{lemma}
\begin{lemma}Let $G$ be a group and  $H$ be a  sugroup of $G$,then $cl_{s_{H}}$ $($ $a$ $)$ is define over $H$ suchthat
           collection of all $cl_{s_{H}}$ $($ $a$ $)$  form a group $($say $K$ $)$ with operation\\ 
$cl_{s_{H}}$ $($ $a$ $)$\textbullet$cl_{s_{H}}$ $($ $b$ $)$=$\{$ c \textasteriskcentered d $/$ $\forall$ $c$ $\in$ $cl_{s_{H}}$ $($ $a$ $)$%and $\forall$ $d$ $\in$ $cl_{s_{H}}$ $($ $b$ $)$ and dot means $.$
is an operation which is equal to operation on $G$ and \textbullet is an operation on $K$ $\}$
iff $a$ \textasteriskcentered $x$=$x$ \textasteriskcentered $a$ $\forall$ $x$ $\in$ $H$ and $a$ $\in$ $G.$
          \end{lemma}
\begin{theorem}
             Let $G$ be a group and $H$ be an faithfull subgroup of $G$ ,$g$ $($ $a$ $)$  define over $H$ suchthat
 collection of all $g$ $($ $a$ $)$ form a group$($say $K$ $)$ and $\exists$ a special normal subgroup $H'$ in $G$ suchthat $K$ $\cong$ $G/H'$.
            \end{theorem}
%\begin{definition}
%\textbf{corrolary}
 %          Let $g$ $($ $1$ $)$ $\in$ $K$ and $g$ $($ $1$ $)$ is an nontrivial derived subgroup of $K$ then $G$ is called \textbf{Camina group}.
   
% and conjugacy classes  $cl(x)$ of $G$ outside the $H'$ coincide with element of $K$ ,suchthat \\ $|$ $cl(x)$ $|$= $|$ $g$ $($ $1$ $)$ $|$
%	    then $G$ is called \textbf{Camina group}.
 %         \end{definition}
%\begin{remark}
% Number of conjugacy classes present inside the  $g$ $($ $1$ $)$ we can classify the \textbf{Camina group}.
          
%\end{remark}
\begin{lemma}
        Let $G$ be an abelian group,$g$ $($ $a$ $)$ define over $G$,then $g$ $($ $1$ $)$ form a normal subgroup of $G$. 
       \end{lemma}
\begin{lemma}
        Let $H$  be a subgroup of $G$,$g$ $($ $1$ $)$ is define over $K$ which is subgroup of $G$ then
    $H$/$($ $H$ $\cap$ $g$ $($ $1$ $)$ $)$ $\cong$ $H$ $g$ $($ $1$ $)$ $/$ $g$ $($ $1$ $)$ .

       \end{lemma}
\begin{lemma}
           Let $N$ be an normal subgroup of $G$,$g$ $($ $1$ $)$ is define over $K$.then\\
	$G$ $/$ $H$ $\cong$ $($ $G$ $/$ $g$ $($ $1$ $)$ $)$ $/$ $($ $H$ $/$ $g$ $($ $1$ $)$ $)$.
          \end{lemma}
%\begin{lemma}
 %       Let $G$ be an abelian group,$g$ $($ $a$ $)$ define over $G$,then $g$ $($ $1$ $)$ form a normal subgroup of $G$. 
  %     \end{lemma}
\begin{lemma}
           Let $G$ be a group and $H$ be an faithfull subgroup of $G$,\\then $O$ $($ $x$\textasteriskcentered$a$\textasteriskcentered$x$ $)$=$O$ $($a$)$
	    if and only if $\forall$ $a$$\in$$H$ ,$\forall$ $x$$\in$$G$ such that $x$\textasteriskcentered$a$=$a$\textasteriskcentered$x$ and
	every element outside the $H$ is self inverse.
          \end{lemma}
\begin{definition}Let $G$ be a Group and $H$ be a nonabelian subgroup of $G$, then\\ $h$ $($ $a$ $)$  is defined as 
$h$ $($ $a$ $)$=$\{$ $x$\textasteriskcentered$a$\textasteriskcentered$x$ $/$  $\forall$ $x$ $\in$ $H$,$a$ $\in$ $G$,$H$ 
be a nonabelian sugroup of $G$ $\}$, then  collection of all $h$ $($ $a$ $)$ is a set $N_{H}$ such that $N_{H}$ form a group then $G$ is 
called as \textbf{Non-R Group}
and $N_{H}$ called as \textbf{Non-selfclass Group} and $h$ $($ $a$ $)$ \textbf{non-self class}.
\begin{remark}
 Notation use to denote non-selfclass of $a$ $\in$ $G$ is $h$ $($ $a$ $)$.
\end{remark}
 % and $a$\textasteriskcentered$x$=$x$\textasteriskcentered$a$ $\forall$ $x$ $\in$ $H$ and $a$ $\in$ $G$ $\}$ 
\end{definition}
\begin{example}
Let $G$ be an non-abelian Group of oder 6 such that element of $G$ is as follow,\\
$G$=$\{$ $E$,$A$,$B$,$C$,$D$,$K$ $\}$ 
then non-self class define over $G$ are  follow,\\
$h$ $($ $A$ $)$=$\{$ $I$\textasteriskcentered$A$\textasteriskcentered$I$,$A$\textasteriskcentered$A$\textasteriskcentered$A$,
$B$\textasteriskcentered$A$\textasteriskcentered$B$,$C$\textasteriskcentered$A$\textasteriskcentered$C$,$D$\textasteriskcentered$A$\textasteriskcentered$D$
,$K$\textasteriskcentered$A$\textasteriskcentered$K$ $\}$\\
$\Rightarrow$ $h$ $($ $A$ $)$=$\{$ $A$,$C$,$K$ $\}$.\\%%%%%%%%%%%%%%%%%%%%%%%%%%%%%%%%%
$h$ $($ $C$ $)$=$\{$ $I$\textasteriskcentered$C$\textasteriskcentered$I$,$A$\textasteriskcentered$C$\textasteriskcentered$A$,
$B$\textasteriskcentered$C$\textasteriskcentered$B$,$C$\textasteriskcentered$C$\textasteriskcentered$C$,$D$\textasteriskcentered$C$\textasteriskcentered$D$
,$K$\textasteriskcentered$C$\textasteriskcentered$K$ $\}$\\
$\Rightarrow$ $h$ $($ $C$ $)$=$\{$ $A$,$C$,$K$ $\}$.\\%%%%%%%%%%%%%%%%%%%%%%%%%%%%%%%%%%%%%%%%
$h$ $($ $K$ $)$=$\{$ $I$\textasteriskcentered$K$\textasteriskcentered$I$,$A$\textasteriskcentered$K$\textasteriskcentered$A$,
$B$\textasteriskcentered$K$\textasteriskcentered$B$,$C$\textasteriskcentered$K$\textasteriskcentered$C$,$D$\textasteriskcentered$K$\textasteriskcentered$D$
,$K$\textasteriskcentered$K$\textasteriskcentered$K$ $\}$\\
$\Rightarrow$ $h$ $($ $K$ $)$=$\{$ $A$,$C$,$K$ $\}$.\\%%%%%%%%%%%%%%%%%%%%%%%%%%%%%%%%%%%%%%%%
$h$ $($ $I$ $)$=$\{$ $I$\textasteriskcentered$I$\textasteriskcentered$I$,$A$\textasteriskcentered$I$\textasteriskcentered$A$,
$B$\textasteriskcentered$I$\textasteriskcentered$B$,$C$\textasteriskcentered$I$\textasteriskcentered$C$,$D$\textasteriskcentered$I$\textasteriskcentered$D$
,$K$\textasteriskcentered$I$\textasteriskcentered$K$ $\}$\\
$\Rightarrow$ $h$ $($ $I$ $)$=$\{$ $I$,$B$,$D$ $\}$.\\
$h$ $($ $B$ $)$=$\{$ $I$\textasteriskcentered$B$\textasteriskcentered$I$,$A$\textasteriskcentered$B$\textasteriskcentered$A$,
$B$\textasteriskcentered$B$\textasteriskcentered$B$,$C$\textasteriskcentered$B$\textasteriskcentered$C$,$D$\textasteriskcentered$B$\textasteriskcentered$D$
,$K$\textasteriskcentered$B$\textasteriskcentered$K$ $\}$\\
$\Rightarrow$ $h$ $($ $B$ $)$=$\{$ $I$,$B$,$D$ $\}$.\\
$h$ $($ $D$ $)$=$\{$ $I$\textasteriskcentered$D$\textasteriskcentered$I$,$A$\textasteriskcentered$D$\textasteriskcentered$A$,
$B$\textasteriskcentered$D$\textasteriskcentered$B$,$C$\textasteriskcentered$D$\textasteriskcentered$C$,$D$\textasteriskcentered$D$\textasteriskcentered$D$
,$K$\textasteriskcentered$D$\textasteriskcentered$K$ $\}$\\
$\Rightarrow$ $h$ $($ $D$ $)$=$\{$ $I$,$B$,$D$ $\}$.\\
clearly, $h$ $($ $A$ $)$=$h$ $($ $K$ $)$=$h$ $($ $C$ $)$and\\ $h$ $($ $I$ $)$=$h$ $($ $B$ $)$=$h$ $($ $D$ $)$.\\
then,collection of all $h$ $($ $A$ $)$ is a set $N_{H}$ such that\\ 
\\$N_{H}$=$\{$ $h$ $($ $A$ $)$,$h$ $($ $I$ $)$ $\}$,\\then $N_{H}$ form a group with respect to operation \textbullet define as follows\\ 
$h$ $($ $A$ $)$\textbullet$g$ $($ $B$ $)$=$\{$ A\textasteriskcentered B $/$ $\forall$ $A$ $\in$ $h$ $($ $A$ $)$ and $\forall$ $B$ $\in$ $h$ $($ $B$ $)$ and dot means \textasteriskcentered
is an operation which is equal to operation on $G$ and \textbullet is an operation on $N_{H}$  $\}$\\
$\Rightarrow$ $N_{H}$ is a \textbf{Non-selfclass Group}\\$\Rightarrow$ $G$ is a \textbf{Non-R Group}
\end{example}
%\begin{remark}

%\end{remark}

\subsection{Rohit conjecture } There exist infinite number of  \textbf{Non-selfclass Group} .

%%%%%%%%%%%%%%%%%%%%%%%%%%%%%%%%%%%%%%%%%%%%%%%%%%%%%%%%%%%%%%%%%%%%%%%%%%%%%%%%%%%%%%%%%%%%%%%%%%%%%%%%%%%%%%%%%%%%%%%%%%%%%%
\section{Conclusion}
\begin{itemize}
\item[1] From self-classes and Conjugacy-selfclass we can study the group structure very clearly.and also 
one can observe that conjugacy class are the special case of Conjugacy-selfclass.
\item[2]one can also observe that size of self-classes and Conjugacy-selfclass depend on subgroup size on which it define.from this we can control 
the size of each class ,so that we can study the property of group.
% \item[\texbf{open problem of nonself class}]
\end{itemize}

\section{Rohit Open Problem on Non-selfclass Group}
\textbf{Description}\\
From definition of non-self class it is clear that $h$ $($ $a$ $)$ not form an equivalence relation  but even 
though there are some group in which all $h$ $($ $a$ $)$ are disjoint,such that collection of all $h$ $($ $a$ $)$ form a group .
The problem is that 
\begin{itemize}
                     \item [1] Does \textbf{Rohit conjecture} is true ?\\
			      if yes ,then which are they?mean what is the oder of every Non-selfclass group,\\
			      if no, then there exist a finite Number of Non-selfclass group,if $S$ be a set of collection of all 
			    Non-selfclass Group,then what is the exact cardinality of set $S$

%does there exit infinite Non-selfclass Group ?
%if yes, which are they, 
%			      if there exit Non-selfclass Group.
                   %  \item[2]  if $G$ be a Non-selfclass Group ,clearly $h$ $($ $a$ $)$ not form an equivalence relation ,
	    %what is the reason behind  collection of all $h$ $($ $a$ $)$ is a set $N_{H}$ such 

                    \end{itemize}
\section{My Future work }
\begin{itemize}
 \item [1] study the property of self-class ,conjugacy self-class,Nonselfclass.
\item[2] study the  \textbf{Rohit conjecture}.
\end{itemize}

%%%%%%%%%%%%%%%%%%%%%%%%%%%%%%%%%%%%%%%%%%%%%%%%%%%%%%%%%%%%%%%%%%%%%%%%%%%%%%%%%%%%%%%%%%%%%%%%%%%%%%%%%%%%%%%%%%%%%%%%%%%%%%%%%%%%%%%%%%%%%%%%%%%%%%%%%%%%%%%%%
\bibliographystyle{amsplain}

\end{document}